\documentclass[reqno,12pt]{amsart}
\usepackage{tabularx}
\usepackage{amsmath, amsthm}
\usepackage{amssymb}
\usepackage{enumitem}
\usepackage{xypic}
\usepackage{mathrsfs}
\usepackage{hyperref}
\usepackage{graphicx}
\usepackage{perpage}
\usepackage{extarrows}
\MakePerPage{footnote}
\usepackage{tikz}
\usepackage{tikz-cd} 
\usetikzlibrary{calc,patterns,angles,quotes}
\usepackage{setspace}
\usepackage{longtable}
\hypersetup{colorlinks,linkcolor={blue},citecolor={blue},urlcolor={red}}  
\allowdisplaybreaks
\usepackage{a4wide}   
\usepackage[final]{pdfpages}

\title[On a link criterion for Lipschitz Normal Embeddings among definable sets]{On a link criterion for Lipschitz Normal Embeddings among definable sets}

\author{Nhan Nguyen}

\address{Basque Center for Applied Mathematics (BCAM),
Alameda de Mazarredo 14, 48009 Bilbao, Bizkaia, Spain}
\email{nnguyen@bcamath.org}
\address{ThangLong Institute of Mathematics and Applied Sciences (TMAS), Nghiem Xuan Yem,  Hoang Mai, Hanoi, Vietnam}
\email{nguyenxuanvietnhan@gmail.com}

\newtheorem{thm}{Theorem}[section]
\newtheorem{lem}[thm]{Lemma}

\newtheorem{defn}[thm]{Definition}

\newtheorem{example}[thm]{Example}
\newtheorem{rem}[thm]{Remark}

\newtheorem{question}[thm]{Question}

\numberwithin{equation}{section}

\newcommand{\bb}{\mathbb}
\newcommand{\al}{\mathcal}

\newcommand{\bff}{\mathbf}

\newcommand{\length}{{\rm length}}
\newcommand{\dist}{{\rm dist}}
\newcommand{\grad}{{\rm grad }}
\newcommand{\id}{{\rm Id }}

\begin{document}
\maketitle
\begin{abstract} 
It is known by a result of Mendes and Sampaio that the Lipschitz normal embedding of a subanalytic germ is fully characterized by the Lipschitz normal embedding of its link. In this note, we show that the result still holds for definable germs in any o-minimal structure on $(\bb R, + , .)$. We also give an example showing that for homomorphisms between MD-homologies induced by the identity map, being isomorphic is not enough to ensure that the given germ is Lipschitz normally embedded. This is a negative answer to the question asked by Bobadilla et al. in their paper about Moderately Discontinuous Homology.
\end{abstract}

\section{Introduction}
Given a connected definable set $X\subset \bb R^n$, one can equip $X$ with two natural metrics: {\it the outer metric ($d_{out}$)} induced by the Euclidean metric of the ambient space $\bb R^n$ and {\it the inner metric ($d_{inn}$)} where the distance between two points in $X$ is defined as the infimum of the lengths of rectifiable curves in $X$ connecting these points. We call $X$ {\it Lipschitz normally embedded} (or {\it LNE} for brevity) if these two metrics are equivalent, i.e., there is a constant $C>0$ such that 
$$\forall x, y \in X, d_{inn} (x,y) \leq C d_{out} (x,y).$$
Any such $C$ is referred to as a {\it LNE constant} for $X$. In order to emphasize the constant $C$ we will call $X$ {\it $C$-LNE}. We say that $X$ is {\it LNE at a point $x_0\in \overline{X}$} (or {\it the germ $(X, x_0)$ is LNE}) if there is a neighborhood $U$ of $x_0$ in $\bb R^n$ such that $X\cap U$ is LNE. It is worth noticing that, in  the o-minimal setting, the definitions of connected and path connected set are equivalent, and every continuous definable curve connecting two points is rectifiable, i.e., it has finite length. 

The notion of Lipschitz Normal Embedding first appeared in a paper of Birbrair and Mostowski \cite{Bir-Mos}, it has since been an active research area and many interesting results were proved, see for example   \cite{Anne2}, \cite{Lev2}, \cite{Bir-Men}, \cite{Bir-Men-Nun}, \cite{Bir-Mos},  \cite{Den-Tib}, \cite{Fer-Sam1}, \cite{Fer-Sam2}, \cite{Ka}, \cite{Anne1}, \cite{Neu-Per-Pic1}, \cite{Neu-Per-Pic2}. Recently, Mendes and Sampaio \cite{Men-Sam} gave a nice criterion for the germ of a subanalytic set to be LNE based on the LNE condition on the link. Namely, they prove that

\begin{thm} [\cite{Men-Sam}, Theorem 1.1]\label{thm_Men-Sam_2} Let $(X,0)$ be a closed connected subanalytic germ in $\bb R^n$. If $(X\setminus \{0\}, 0)$ is connected, then the following statements are equivalent: 

(1) $(X, 0)$ is LNE;

(2) $(X, 0)$ is LLNE.
\end{thm}

If $(X\setminus \{0\}, 0)$ has many connected components, the result is generalized as follows: 

\begin{thm}[\cite{Men-Sam}, Theorem 4.1] \label{thm_Men-Sam}
Let $(X,0)$ be a closed connected subanalytic germ in $\bb R^n$. Suppose that $\{X_i\}_{i\in I}$ are connected components of $(X\setminus \{0\},0)$. Then, the following statements are equivalent:

(1) $(X,0)$ is LNE.

(2)  $(\overline{X}_i, 0)$ is LNE for every $i\in I$ and there is  $C>0$ such that $$\min_{j, k \in I, j\neq k}\{\dist(L_{r}(X_j), L_{r}(X_k))\} \geq C r$$ for $r$ sufficiently small.

(3) $(\overline{X}_i,0)$ is LLNE for every $i\in I$ and there is  $C>0$ such that $$\min_{j, k \in I, j\neq k}\{\dist(L_{r}(X_j), L_{r}(X_k))\} \geq C r$$ for $r$ sufficiently small.

Here, $L_{r}(X_j)$ and $L_{r}(X_k)$ are the $r$-links of $X_j$ and $X_k$ respectively, and $\dist(., .)$ is the usual outer distance. 
\end{thm}

Let us recall the notions of {\it $r$-link} and {\it link Lipschitz normal embedding (LLNE)} of a definable germ. Let $(X, 0)$ be a definable germ in $\bb R^n$. Given $r >0$, the {\it $r$-link} of $X$ is the set $L_{r}(X) := X \cap \bff S^{n-1}_r$ where $ \bff S^{n-1}_r$ denotes the $(n-1)$-dimensional sphere of radius $r$ centered at the origin. By the locally conic structure (see for example \cite{Coste}, Theorem 4.10) the definable topological type of the $r$-link is invariant when $r$ is small enough. 
One calls an $r$-link of $(X, 0)$ for $r$ sufficiently small {\it the link} of the germ $(X,0)$ (or {\it the link associated to the Euclidean norm}).

 In general, given a continuous definable function germ $\rho: (\bb R^n, 0) \to \bb R_{\geq 0}$ such that $\rho^{-1}(0) = \{0\}$, for $r>0$, {\it the $r$-link of $(X, 0)$ associated to $\rho$} is the set $L_{\rho, r} (X) : = \rho^{-1}(r) \cap X$. It is known that for $r$ sufficiently small, the definable topological type of $L_{\rho, r} (X)$ is uniquely determined, it depends neither on $r$ nor $\rho$ (\cite{Kur}, Proposition 1).  We call such an $r$-link for $r$ small enough {\it the general link} or {\it the link associated to $\rho$} in case we want to emphasize the function $\rho$. 
 
 Given a definable germ $\rho: (X, 0) \to \bb R_{\geq 0}$ which allows a continuously definable extension to a definable function germ $\rho': (\bb R^n, 0) \to \bb R_{\geq 0}$ such that $\rho'^{-1}(0) = \{0\}$,  then we can define {\it the link of $(X, 0)$ associated to $\rho$} to be the set $L_{\rho, r}:= \rho^{-1}(r)\cap X$ for $r$ sufficiently small.
 
 Suppose that  $(X\setminus \{0\}, 0)$ is connected. For a continuous definable function $\rho: (X, 0) \to \bb R_{\geq 0}$ given as above, we say that $(X, 0)$ is {\it LLNE with respect to $\rho$} if there are $r_0 >0$ and $C >0$ such that for all $0< r< r_0$, $L_{\rho, r}(X)$ is LNE with an LNE constant bounded above by $C$. We call such a constant $C$ an {\it LLNE constant} (w.r.t $\rho$) for $(X,0)$. To emphasize the constant $C$, we call $X$ {\it $C$-LLNE with respect to $\rho$}. We simply call $(X, 0)$ {\it LLNE} if it is LLNE with respect to the Euclidean norm.  
 
The arguments of  Mendes--Sampaio in the proof of Theorem \ref{thm_Men-Sam_2} rely first on a  result of Birbrair--Mendes \cite{Bir-Men}, which says that for semialgebraic germs, being LNE is equivalent to the condition that for any pair of arcs parametrized by distance to the origin, the inner and outer contacts are the same, and then on the following result due to Valette:
\begin{thm}[\cite{Valette3}, Theorem 4.4.8; \cite{Valette2}, Corollary 2.2] \label{thm_Valette} Let $(X,0)$ be a subanalytic germ in $\bb R^n$.  Let $\rho: (X,0) \to (\bb R,0)$ be  a subanalytic radius function. Then there is a bi-Lipschitz subanalytic homeomorphism $h: (X,0) \to (X,0)$ such that  $\|h(x)\| = \rho(x)$.
\end{thm}
Here, a function germ $\rho: (X,0) \to (\bb R,0)$ is called {\it a radius function} if it is Lipschitz and $\rho(x) \sim \|x\|$. 

In fact, Mendes--Sampaio showed further that their result is still true for the link associated to a subanalytic norm.  One may check that the same method can be applied to definable germs in polynomially bounded o-minimal structures. However, improvement is needed to make it work on an arbitrary o-minimal structure. The reason for this is two-folds:  first the concept of order between two definable arcs as in Birbrair--Mendes' result does not make good sense and second the result of Valette is no longer valid in this context. Theorem \ref{thm_Valette} implies that, up to a bi-Lipschitz subanalytic homeomorphism, the link of associated to a subanalytic radius function is uniquely determined. Consequently, given a radius function, for any $0< r_1< r_2 < r$ for $r$ small enough, the sets $L_{\rho, r_1}(X)$ and $L_{\rho, r_2}(X)$ are subanalytically bi-Lipschitz equivalent (of course, the Lipschitz constant depend on $(r_1, r_2)$). This property holds only for definable germs in polynomially bounded o-minimal structures, it fails in the exponential o-minimal setting. Consider Parusi\'nski's example  $X = \{(x, y, t) \in \bb R^3: |y| = x^{t + 1},0\leq  x \leq t\}$ and $\rho_1(.) = \|.\|_{\infty}$. It is easy to check the Lipschitz types of $r$-links of $(X,0)$ associated to $\rho(x) = \|x\|_\infty$ vary continuously.

In this paper, using different techniques, we show that Theorem \ref{thm_Men-Sam}, hence Theorem \ref{thm_Men-Sam_2}, holds for definable germs in any o-minimal structure. Furthermore, the results remain valid if we replace the link defined by the Euclidean norm with the link associated to a  {\it definable radius function}. More precisely,  for an arbitrary fixed o-minimal structure, we prove that

\begin{thm}[Main Theorem]\label{thm_main}
Let $(X,0)$ be a connected definable germ in $\bb R^n$.  Let $\rho: (X,0) \to (\bb R,0)$ be a definable radius function. Suppose that $\{X_i\}_{i\in I}$ are connected components of $(X\setminus \{0\},0)$. Then, the following statements are equivalent 

(i) $(X,0)$ is LNE.

(ii) $(X_i, 0)$ is LNE for every $i \in I$ and there is  $C>0$ such that $$\min_{j,k \in I, j\neq k} \{\dist(L_{\rho, r}( X_j), L_{\rho, r}( X_k)\} \geq C r$$ for $r$ sufficiently small.

(iii) $(X_i,0)$ is LLNE with respect to $\rho$ for every $i \in I$ and there is  $C>0$ such that $$\min_{j,k \in I, j\neq k} \{\dist(L_{\rho, r}( X_j), L_{\rho, r}( X_j))\} \geq C r$$ for $r$ sufficiently small.
\end{thm}

A direct consequence of Theorem \ref{thm_main} is that:
  \begin{thm}\label{thm_main_2}
Let $(X,0)$ be a connected definable germ in $\bb R^n$.  Let $\rho: (X,0) \to (\bb R,0)$ be a definable radius function. Suppose that $(X\setminus\{0\}, 0)$ is connected. Then, $(X,0)$ is LNE if and only if it is LLNE with respect to $\rho$.
\end{thm}

Remark that every germ of a definable norm is a definable radius function. The converse is not true, for example, the function $\rho + \rho^2$, where $\rho$ is a definable norm, is a radius function but no longer a norm. In Theorem \ref{thm_main}, $X$ need not to be closed.

In \cite{Bobadilla} Bobadilla--Heinze--Pereira--Sampaio introduced a homology theory called Moderately Discontinuous homology (MD-homology) in order to capture the singular homology of the link of a given subanalytic germ $(X,0)$ after collapsing it at a certain speed.  The identity map $\id_{(X,0)}: (X, 0, d_{inn}) \to (X,0, d_{out})$  induces homomorphisms between groups of MD-homologies of $(X,x_0)$ for all $x_0 \in (X,0)$. It is easy to check that if $(X,0)$ is LNE then these homomorphisms are actually isomorphisms. It is  asked in \cite{Bobadilla} if the converse holds, i.e., suppose these homomorphisms are isomorphisms then is it true that $(X,0)$ is LNE. In Section \ref{section3}, as an application of Theorem \ref{thm_Men-Sam} (also Theorem \ref{thm_main_2}) we give a simple example showing that in general the answer is negative. 

Throughout the paper, we assume that the reader is familiar with the notion of o-minimal structures on $\bb R$. By ``definable" we mean definable in  a given o-minimal structure. We will use Curve Selection (\cite{Coste}, Theorem 3.2), Definable Choice (\cite{Coste}, Theorem 3.1) without reciting the references.  We refer the reader to  \cite{Coste}, \cite{Dries} for the details about the theory of o-minimal structures. 

We denote by $\bff B^{n}_r$, $\bff S^{n-1}_r$ respectively the $n$-dimensional closed ball and the $(n-1)$-dimensional sphere in $\bb R^n$ of radius $r$ centered at $0$. Let $X\subset \bb R^n$.  We denote by $\overline{X}$ the closure of $X$ in $\bb R^n$. Given non-negative functions $f, g: X \to \bb R$, we write $f\lesssim g$ (or $g\gtrsim f$) if there is $C>0$ such that $f(x) \leq C g(x)$ for every $x \in X$. Such a  constant $C$ is called {\it a constant for the relation $\lesssim$}.  We write $f\sim g$ if $f\lesssim g$ and $g\lesssim f$, i.e., there are constants $C_1, C_2>0$ such that 
$C_1 f(x) \leq g(x) \leq C_2 f(x)$ for every $x \in X$.  Given two continuous function germs $f, g : (X,0) \to \bb R$ such that $g(0)\neq 0$, we write $f\ll g$ if $\lim_{x\to 0} f(x)/g(x) \to 0$.  Let $Y$ be another subset of $\bb R^n$. We denote by $\dist(X,Y)$ the usual outer distance between $X$ and $Y$. We use $\|.\|$  and $\|.\|_{\infty}$ respectively for the Euclidean norm and the maximum norm in $\bb R^n$. We also use the notation $d_X$ for the inner metric on $X$.

\section{Proof of the main theorem}\label{section2}
This section is devoted to the proof of Theorem \ref{thm_main}. The idea of the proof is the following.  We first extend $\rho$ to a definable radius function defined on the whole of $(\bb R^n, 0)$, we use the same notation $\rho$ for this extension.  Note that the map germ  $\varphi: (\bb R^n, 0) \to (\bb R^n, 0)$ defined by $\varphi(x):= \frac{\rho(x)}{\|x\|_{\infty}} x$ if $x \neq 0$ and by $\varphi(0) := 0$  brings $L_{\rho, r}(X)$ to $ L_{\rho', r}(\varphi(X))$ where $\rho'(x) : = \|x\|_\infty$. By showing that $\varphi$ is the germ of a bi-Lipschitz homeormophism (see Lemma \ref{lem_bi-Lipschitz}), we can reduce our problem to proving the theorem for the link with respect to $\rho(x) = \|x\|_\infty$. 

In Lemma \ref{lem_conic2}, we construct a vector field whose flow provides a ``shortest" way to go from a point in $L_{\rho, r}(X) $ to $L_{\rho, r'}(X)$ where $r\neq r'$. This is sufficient for us  to prove that $(iii) \Rightarrow (ii) \Rightarrow (i)$. 

To show that  $(i) \Rightarrow (iii)$, we just need to consider the case that  $(X\setminus\{0\},0)$ has only one connected component (the proof for the case of many connected components follows easily by Curve Selection). As in the proof of Mendes--Sampaio, we separate $(X, 0)$ into finitely many Lipschitz cells which are proved to be LLNE with respect to $\rho$ (see Lemma \ref{lem_LLNE}) (to make a proof simpler, we ask each of these Lipschitz cells after a linear change of coordinates is contained in a special area denoted by $\al C_n$). The vector field in Lemma \ref{lem_conic2} can be chosen to be compatible with these Lipschitz cells. This allows us to construct a special path connecting two given points $p, q \in L_{\rho, r}(X)$ whose length is equivalent to $d_X(p, q)$. The result then follows.  

We now start with some preliminaries. 

\begin{defn}\rm Let $C\subset \bb R^n$ be a definable set. 
 
 (1) We call $C$ a {\it standard Lipschitz cell  with constant $M$} if 
 
     $n=1$: $C$ is a point or an open interval, 
     
     $n>1$: $C$ has one of the following forms
     \begin{itemize}
         \item $\Gamma_{\xi}: =\{(x,y)\in B\times \bb R:  x \in B, y = \xi(x)\}$
         \item $(\xi_1, \xi_2): = \{(x,y)\in B\times \bb R:  x \in B, \xi_1(x)< y < \xi_2(x)\}$
         \item $(-\infty, \xi): = \{(x,y)\in B\times \bb R:  x \in B, y < \xi(x)\}$
         \item $(\xi, +\infty): = \{(x,y)\in B\times \bb R:  x \in B, \xi(x)< y\}$
     \end{itemize}
      where $B\subset \bb R^{n-1}$ is a standard Lipschitz cell with constant $M$ and $\xi, \xi_1, \xi_2$ are $C^1$-definable functions over $B$ with the first derivatives bounded above by $M$. We call $B$ the {\it basis of $C$}. If we do not care much about the constant $M$ we simply call $C$ a {\it standard Lipschitz cell}.

   (2) We call $C\subset \bb R^n$ an {\it Lipschitz cell (with constant $M$)}  if there is an orthogonal change of coordinates $\phi: \bb R^n \to \bb R^n$ such that $\phi(C)$ is a standard Lipschitz cell (with constant $M$). 
 \end{defn}
 \begin{rem}\label{rem_L_cell} \rm  Let $C\subset \bb R^n$ be a Lipschitz cell with constant $M$. The following results follow directly from the definition:
 \begin{enumerate}
     \item $C$ and $\overline{C}$ are LNE with LNE constants depending only on $M$. 
     \item  For any two points $x, y$ in $\overline{C}$ there is a continuous definable curve $\gamma: [0,1]\to \overline{C}$ such that $\gamma(0) = x$, $\gamma(1) = y$, $\gamma((0,1)) \subset C$ and 
 $$\length(\gamma) \sim d_{\overline{C}}(x,y)\sim  \|x-y\|,$$
 \end{enumerate}
 where  constants of the equivalence relation $\sim$ depend only on $M$.
 \end{rem}
 The following result is due to Kurdyka \cite{Kur} in the subanalytic category, his proof also works for o-minimal structures.  One may find a parameterized version of the theorem in \cite{Kur-Par}, Proposition 1.4.
 
 \begin{thm}\label{thm_ku_par} Let $\al X= \{X_1, \ldots, X_k\}$ be a finite family of definable sets in $\bb R^n$. Then, there is a finite partition of $\bb R^n$ into Lipschitz cells compatible with $X_i$'s, i.e., each $X_i$, $i = 1, \ldots,k$ is a union of some elements of the partition.
 \end{thm}

We now define the special area $\al C_n$  mentioned at the beginning of the section. For $n \in \bb N$, we set
\begin{equation}\label{fm_C_n}\al C_n := \{(x_1, \ldots, x_n)\in \bb R^n: \|x\|_{\infty} = x_1\}.
\end{equation}
Note that for $i = 1, \ldots, n$, $\pi_i(\al C_n) = \al C_i$  where $\pi_i: \bb R^n \to \bb R^i$ is the projection onto the first $i$ coordinates. 

 \begin{lem}\label{lem_LLNE} Let $C$ be a Lipschitz cell in $\al C_n$ such that $\dim C \geq 1$ and $0\in  \overline{C}$. Then, $(C,0)$ and $(\overline{C},0)$ are LLNE  with respect to $\rho (x) := \|x\|_{\infty}$.
 \end{lem}
 \begin{proof} Since $\dim C \geq 1$ and $C\subset \al C_n$,  $\pi_1(C)$ is an interval $(0, a)$ for some $a>0$. For $0< r< a$, we set $C(r) := \rho^{-1}(r) \cap C$ and $ \overline{C}(r) := \rho^{-1}(r)\cap \overline{C}$. Since $C \subset \al C_n$ is a Lipschitz cell with a constant $M$, $C(r)$ is  also a Lipschitz cell with constant $M$. Moreover, $\overline{C}(r)$ is the  closure of $C(r)$,   by Remark \ref{rem_L_cell}, part (2), if $\overline{C}(r)$ is LNE with LNE constant independent of $r$ so is  $C(r)$. Consequently, if  $(\overline{C}, 0)$ is LLNE then $(C, 0)$ is also LLNE.
 Therefore, to prove the lemma it suffices to show that  $(\overline{C}, 0)$ is LLNE. 
 
 The proof is by induction on $n$.  For $n = 1$, the result is trivial.
 
 Assume that $n \geq 2$. Let $B$ be the basis of $C$. Then, $B\subset \al C_{n-1}$, so  by inductive assumption, $(\overline{B}, 0)$ is LLNE with respect to $\rho$. There are two possibilities for $C$, it is either a graph or a band over $B$.
 
 \underline{\it Case 1:}  $C = \Gamma_\xi$ (graph). Since $\xi$ is Lipschitz, it is possible to extend it to a Lipschitz function $\overline{\xi}$  over $\overline{B}$. Since $\overline{C}(r)$ is the graph of the restriction of $\overline{\xi}$ to $\overline{B}(r) := \overline{B} \cap \rho^{-1}(r)$ and  $\overline{B}(r)$ is LNE with an LNE constant independent of $r$, it follows that  $\overline{C}(r)$ is also LNE with an LNE constant independent of $r$.  
 
 \underline{\it Case 2:} $C = (\xi_1, \xi_2)$ (band). Extend  $\xi_1, \xi_2$ to Lipschitz functions $\overline{\xi}_1, \overline{\xi}_2$ on $\overline{B}$. It is clear that $ \overline{C} = \{(x,y)\in \overline{B}\times \bb R : \overline{\xi}_1(x) \leq y \leq \overline{\xi}_2(x)\}$ and 
 $$\overline{C}(r) = [\overline{\xi}_1|_{\overline{B}(r)}, \overline{\xi}_2|_{\overline{B}(r)}]: = \{(x, y) \in \overline{B}(r) \times \bb R:  \overline{\xi}_1(x) \leq y \leq \overline{\xi}_2(x)\}.$$
 
Let $z$ and $z'$ be two points in $\overline{C}(r)$. There are finite points $\{z = z_0, z_1, \ldots, z_k = z'\}$ contained in the intersection of the segment $[z, z']$ and the set $\Gamma_{\overline{\xi}_1|_{\overline{B}(r)}} \cup \Gamma_{\overline{\xi}_2|_{\overline{B}(r)}}$ such that for each $i< k$, the open interval $(z_i, z_{i+1})$ is contained either in $\overline{C}(r)$ or in $(\al C_n\cap \{x_1 = r\}) \setminus \overline{C}(r)$. If $(z_i, z_{i+1})\subset \overline{C}(r)$ then $\|z_i - z_{i+1}\| = d_{\overline{C}(r)}(z_i, z_{i+1})$. Recall that $d_{\overline{C}(r)}$ is the inner distance on $\overline{C}(r)$. If $(z_i, z_{i+1})\subset (\al C_n\cap \{x_1 = r\}) \setminus \overline{C}(r)$ then $z_i$ and $z_{i+1}$ are both in the same graph of the restriction $\overline{\xi}_j|_{\overline{B}(r)}$, $j\in \{ 1,2\}$. Since $\overline{B}(r)$ is LNE with a constant independent of $r$ and $\overline{\xi}_j$ is Lipschitz, the graph of  $\overline{\xi}_j|_{\overline{B}(r)}$ is LNE with an LNE constant independent of $r$.  This implies that 
$$\|z_i - z_{i+1}\| \sim d_{\Gamma_{\overline{\xi}_j|_{\overline{B}(r)}}}(z_i , z_{i+1})  \geq d_{\overline{C}(r)}(z_i , z_{i+1}).$$ 
Therefore, 
 $$\|z - z'\|  = \sum_{i = 0}^{k-1}\|z_i - z_i\| \gtrsim \sum_{i = 0}^{k-1}d_{\overline{C}(r)} (z_i , z_{i+1}) \geq d_{\overline{C}(r)} (z , z').$$
 Hence, $\overline{C}(r)$ is LNE with an LNE constant independent of $r$. Consequently, $(\overline{C}, 0)$ is LLNE.
 \end{proof}

\begin{lem}\label{lem_conic2}
Let $\al X = \{X_i\}_{i\in I}$ be a finite collection of definable germs at $0$. Let $\rho: (\bb R^n,0) \to (\bb R,0)$ be a radius function. Then, there are $r_0>0$, $C>0$, a definable $C^2$ stratification $\Sigma$ of $ X: = \bb R^n\setminus\{0\}$ compatible with $\{X_i\setminus{0}\}_{i\in I}$ and a continuous integrable stratified vector field $\xi$ on $\Sigma$ such that the flow $\Phi: U\subset  X \times \bb R \to X $ generated by $\xi$ has the following properties for all $x\in X \cap \rho^{-1}((0, r_0])$:

(1) $\Phi$ preserves the strata of $\Sigma$, i.e., for each $x \in S\in \Sigma$, $\phi_x(s) \in S$ for every $s\in U_x: =  U \cap (\{x\}\times \bb R)$;

(2) $U_x$ contains the interval  $[0, \rho(x))$;
   
 For every  $0\leq  s <  \rho(x)$, we have
    
(3) $ \Phi_x (s) \in X \cap  \rho^{-1}(r-s)$;
     
(4) $\length(\Phi_x ([0, s]) \leq Cs$.
  
\end{lem}

\begin{proof}
Let $\Sigma$ be a Whitney $(b)$-regular stratification of $X$ compatible with  $\{X_i\setminus \{0\}\}_{i\in I}$ such that the restriction of $\rho$ to each stratum of $\Sigma$ is of class $C^2$. The existence of  Whitney stratification for definable sets is proved in \cite{Loi2} (see also \cite{Loi1}, \cite{Ngu-Tri-Tro}). For $x\neq 0$,  set $v(x) := \frac{x}{\|x\|}$. For $x \in S \in \Sigma$, set $w(x) := P_x (v(x))$ where $P_x: \bb R^n \to T_x S$ is the orthogonal projection from $\bb R^n$ to the tangent space to $S$ at the point $x$. 

We claim that 
 
 (i) Given $\varepsilon>0$ there is $R>0$ such that for every  $x$ in  $X \cap\bff B^{n}_R $ we have
$$\| w(x) - v(x)\| < \varepsilon.$$

(ii)  There are $\varepsilon'$ and $R >0$  such that for every $ x\in S\cap \bff B^{n}_R, S\in \Sigma$, we have 
$$ D_x(\rho|_S) (w(x)) > \varepsilon'$$
where $D_x(\rho|_S)$ denotes the tangent map at $x$ of the restriction $\rho|_S$. 

We start with a proof of (i). Assume that  (i) were not true. By Curve Selection, there are  $\varepsilon>0$, a stratum $S\in \Sigma$ and a $C^1$ definable curve $\gamma: [0,\delta) \to \bb R^n$ with $\gamma(0) =0$, $\gamma((0, \delta)) \subset S$ such that for every $t\in (0, \delta)$:
$$\|w(\gamma(t)) - v(\gamma(t)\| \geq \varepsilon \text{ for every } t \text{ in } (0, \delta).$$
Since $\gamma$ is a $C^1$ curve through the origin, the angle between $v(\gamma(t)) = \frac{\gamma(t)}{\|\gamma(t)\|}$ and the tangent line to $\gamma$ at $\gamma(t)$ tends to $0$ as $t$ tends to $0$. Therefore, the angle between $v(\gamma(t))$ and the tangent space $T_{\gamma(t)} S$ tends to $0$ as well, which gives a contradiction.

Now we prove (ii).  We first prove that   
\begin{equation}\label{equ2.3.1}
    |D_x(\rho|_S) (w(x))| > \varepsilon' \text{ for every } x \text{ in }  S\cap \bff B^{n}_R.
\end{equation}
Note that since $\rho$ is Lipschitz, there is $M>0$ such that $\|D_x(\rho|_S)\|< M$  for every $x\in S$ and for every $S\in \Sigma$. This implies that for any $C^1$-definable curve $\gamma:  [0, \delta) \to \bb R^n$ such that $\gamma(0) = 0$ and $\gamma((0, \delta)) \subset S$ the limit   $\lim_{t \to 0} D_{\gamma(t)} (\rho|_S)$  always exists.

We assume on the contrary that (ii) fails. By Curve Selection there are  $S\in \Sigma$  and a  $C^1$-definable curve $\alpha: [0, \delta) \to \bb R^n$ with $\alpha(0) = 0$ and $\alpha((0, \delta))\subset S$ such that $D_{\alpha(t)}(\rho|_S) (w(\alpha(t))$ tends to $0$ as $t \to 0$. Reparametrizing $\alpha$ if necessary, we may assume $\|\alpha (t)\| \sim t$. Then $\alpha(t) = a t + o(t)$ where $a \in \bb R^n, a \neq 0$. Note that $\lim_{t\to 0}v(\alpha(t))= \lim_{t\to 0}w(\alpha(t)) = \lim_{t\to 0} \frac{\alpha'(t)}{\|\alpha'(t)\|},$ where $\alpha'(t)$ denotes the derivative of $\alpha$ at $t$. Therefore,  
\begin{align*}
   0= \lim_{t\to 0} D_{\alpha(t)}(\rho|_S) (w(\alpha(t)) & = \lim_{t\to 0} D_{\alpha(t)}(\rho|_S) \left(\frac{\alpha'(t)}{\|\alpha'(t)\|}\right) \\
    & = \lim_{t\to 0} \frac{(\rho \circ \alpha)'(t)}{\|\alpha'(t)\|}.
\end{align*}
Since $\lim_{t\to 0}\|\alpha'(t)\| = \|a\| \neq 0$,  $\lim_{t\to 0}  (\rho \circ \alpha)'(t) = 0$. This implies that $(\rho\circ \alpha)(t) \ll t$, which contradicts the fact that  $(\rho\circ\alpha)(t) \sim \|\alpha (t)\| \sim t$. Hence, (\ref{equ2.3.1}) is proved.

Note that since $\rho(x) \sim \|x\|$, for any curve $\beta: [0, \delta) \to \bb R^n$ $\beta(0) = 0$, $\beta((0, \delta)) \subset S \in \Sigma$, $(\rho\circ \beta) (t)$ is strictly increasing for $t$ near $0$.  Thus, $D_{\beta(t)}(\rho|_S) (w(\beta(t))>0$ for $t$ near $0$.  Combining with (\ref{equ2.3.1}),  we have $D_x(\rho|_S) (w(x))>0$ for every  $ x\in S\cap \bff B^{n}_R$, so (ii) is proved.

Now we are in a position to prove the lemma. It is shown in the proof of Lemma 3.2 in \cite{ng-va} that for any $\varepsilon >0$, there are $R>0$ and a continuous integrable stratified vector field $\mu$ on $\Sigma \cap \bff B_R^n$ such that $\|\mu(x) - v(x)\|< \varepsilon$ (in fact, the Whitney condition $(b)$ is used in this step). By (i) in the claim, shrinking $R$ if necessary, we may assume that $\|\mu(x) - w(x)\| < 2\varepsilon$. By (ii), $D_x(\rho|_S) (w(x)) > \varepsilon'$ for some $\varepsilon'>0$. Taking $\varepsilon$ sufficiently small, we have $D_x(\rho|_S) (\mu(x)) > \varepsilon'/2$.

For $x \in S\in \Sigma$, set 
$$ \xi(x) := \frac{-\mu(x)}{D_x(\rho|_S) (\mu(x))}.$$ 
For $\varepsilon$ small enough, $\|\mu(x)\|\sim \|w(x)\|\sim \|v(x)\|$. Since $\|v(x)\|=1$, there is $C>0$ such that $\|\xi\|<C$.  Furthermore, since $\mu$ is integrable so is $\xi$. 

Let $\Phi$ denote the flow generated by $\xi$. Take $r_0 >0$ small such that $\rho^{-1}((0, r_0)) \subset \bff B^n_R$ (this is possible because $\rho(x)\sim \|x\|$). We now check that $\Phi_x$ satisfies the conditions (1)--(4). The condition (1) is clear since $\xi$ is a stratified vector field on $\Sigma$.  Note that for every $x \in S \cap \bff B^n_R $, $S\in \Sigma$, we have 
\begin{equation}\label{equ_vector_field}
    D_x(\rho|_S)(\xi(x)) = -1.
\end{equation} 
This implies that on $X\cap \bff B_R^n$,  $\xi$ is a lifting of the constant vector field $\eta(z) = -1$ defined on the image of the set $X  \cap \bff B_R^n$ under the map $\rho$. The domain of $\Phi_x$ then coincides with the domain of the flow generated by the vector field $\eta$ at the point $z = \rho(x)$ which obviously contains the interval $[0, \rho(x))$.  Thus, (2) is proved.

Again, by (\ref{equ_vector_field}),
for $x \in S \cap \rho^{-1}(r)$ ($r\leq r_0$) and $s<r$, one has $\Phi_x(s) \in S \cap \bff \rho^{-1}(r-s)$. Condition (3) then follows.  Finally,  $$\length(\Phi_x  ([0,s])) = \int_0^s |\xi(\Phi_x(s)|ds \leq C \int_0^s ds = Cs,$$ and hence (4) is satisfied.
\end{proof}

The following results are proved in \cite{Valette3} for subanalytic germs, the same arguments work also in the o-minimal setting. 

\begin{lem}[\cite{Valette3}, Lemma 4.4.4]\label{lem_val_1} Let $(X, 0)$ be definable germ in $\bb R^n$ and let $\rho: (X, 0) \to (\bb R, 0)$ be a definable radius function. Let $x, y : (0, \delta)\to X$  be two definable curves such that $\rho(x(t)) = \rho(y(t)) = t$ for every $0< t <\delta$. Then, if $x$ and $y$ have the same tangent cone at the origin then the lines $\overline{0x(t)}$ and $\overline{x(t)y(t)}$ are of different limits when $t$ tends to $0$. 

\end{lem}
\begin{lem}[\cite{Valette3}, Lemma 4.4.6]\label{lem_val_2} Let $f: (X, 0) \to (Y, 0)$ be a definable map germ.  Let $\alpha: (Y, 0) \to \bb R$ be a definable radius function such that $\rho(x): = \alpha(f(x))$ defines a radius function on $(X, 0)$. If the function $f|_{\rho^{-1}(t)\cap X}$ is $L$-bi-Lipschitz  for every $t>0$ small with $L$ indepedent of $t$, then $f$ is the germ of a bi-Lipschitz map. 
\end{lem}

\begin{lem}\label{lem_bi-Lipschitz}
Let $\rho: (\bb R^n,0) \to ( \bb R,0)$ be a definable radius function. Let $\varphi: (\bb R^n,0) \to (\bb R^n,0)$ be defined by $\varphi(x):= \frac{\rho(x)}{\|x\|_\infty} x $ if $x \neq 0$ and $\varphi(0) = 0$. Then, $\varphi$ is the germ of a bi-Lipschitz definable homeomorphism.
\end{lem}
\begin{proof}   First we show that $\varphi$ is a bijection. 
Note that $x$  and $\varphi(x)$ lie in the same half line starting at the origin.  To prove the bijectivity of  $\varphi$ it suffices to show that there is $R>0$ such that for any $l\in \al L$ the restriction of $\varphi$ to $l\cap (\bff B^n_R\setminus \{0\})$ is injective. Here  $\al L$ denotes the set of all half lines in $\bb R^n$ starting at the origin.

{\it Claim:} There is $R>0$ such that for every $l \in \al L$ the restriction 
$\rho|_{l\cap \bff B_R^n}: (l\cap \bff B_R^n,0)   \to (l\cap \bff B_R^n,0)$ is bijective. 

On the contrary we assume that the claim  were not true.  By Curve Selection, there are two $C^1$ definable curves $\gamma_1, \gamma_2: [0, \varepsilon)\to \bb R^n$ such that 
 
 (1) $\gamma_1(0) = \gamma_2(0) = 0$, $\gamma_1((0, \varepsilon))$ and $\gamma_2((0,\varepsilon))$ are in $\bb R^n\setminus \{0\}$
 
 (2) $\gamma_1(t)\neq \gamma_2(t)$ and they are both contained in some $l_t \in \al L$ for every $t\in (0, \varepsilon)$, 
 
 (3) $\rho(\gamma_1(t)) = \rho(\gamma_2(t))$ for every $t\in (0, \varepsilon)$.

It follows from (1) and (2) that the lines $\overline{0 \gamma_1(t)}$ and $\overline{\gamma_1(t) \gamma_2(t)}$ have the same limit in the Grassmannian when $t$ tends to $0$. From (3), by reparametrizing,  we may assume that $\rho(\gamma_1(t)) = \rho(\gamma_2(t)) = t$ for every $t\in (0, \varepsilon)$. By Lemma \ref{lem_val_1},  $\overline{0 \gamma_1(t)}$ and $\overline{\gamma_1(t) \gamma_2(t)}$ have different limits, which is a contradiction, so  the claim is proved.

We now take $R>0$ as in the claim.  Let $x$ and $x'$ be in $l\cap \bff (\bff B^n_R\setminus \{0\})$. There is $t\neq 0$ such that $x' = tx$. We have 
\begin{align*}
    \varphi(x) = \varphi (x') &\Leftrightarrow \frac{\rho(x)}{\|x\|_\infty} x =  \frac{\rho(tx)}{\|tx\|_\infty} tx\\
    &\Leftrightarrow\rho(x) = \rho(tx) = \rho (x') \Leftrightarrow x  = x' \text{ (since } \rho|_{l\cap \bff B^n_R} \text{ is injective)} .
\end{align*}
Therefore, $\varphi|_{l\cap \bff B^n_R\setminus \{0\})}$ is injective.

To prove $\varphi$ is bi-Lipschitz on $\bff B^n_R$, by Lemma \ref{lem_val_2},  it is enough to show $\varphi|_{\rho = r}$ is $L$-bi-Lipschitz for every $r$ small with $L$ independent of $r$.

 Let $x , x' \in \bb R^n$ such that $\rho(x) = \rho(x') = r$. Since $\rho$ is a radius function, $\|x\| \sim \|x'\| \sim r$, and hence $\|x \|_\infty \sim \|x'\|_\infty \sim r$. Then, 
\begin{align*}
   \| \varphi(x) - \varphi(x')\| & = r \bigg\| \frac{x}{\|x\|_\infty} - \frac{x'}{\|x'\|_\infty}\bigg\|\\
   & \leq \frac{r}{\|x\|_\infty} \|x -x'\| + \frac{ r \|x'\|}{\|x\|_\infty \|x'\|_\infty} | \|x'\|_\infty - \|x\|_\infty| \\
    & \lesssim \|x - x'\| +  \|x-x'\|_\infty \\
    &\lesssim \|x - x'\|.
\end{align*}
Therefore, $\varphi|_{\rho = r}$ is $L$-Lipschitz with some constant $L$ independent of $r$.

Now we show that $\|\varphi(x) - \varphi(x')\| \gtrsim \|x - x'\|$. Assume on the contrary that there are two definable curves $x, x': [0, \delta) \to \bb R^n$ such that $x(0) = x'(0) = 0$, $\rho(x(t)) = \rho(x'(t))$ such that 
\begin{equation}\label{equ_kk}
    \lim_{t \to 0} \frac{\| \varphi(x(t)) - \varphi(x'(t))\|}{\|x(t) - x'(t)\|} \to 0.
\end{equation}
The proof is split into two cases: 

{\it Case 1:} $x$ and $x'$ have different tangent cones at the origin. It is easy to see that $\|x(t) - x'(t)\|\sim \|x\| \sim r $ and $\|\varphi(x(t)) - \varphi(x(t'))\| \sim \|\varphi(x(t))\| \sim r$, which contradicts (\ref{equ_kk}).

{\it Case 2:} $x$ and $x'$ have the same tangent cone at the origin. Let $l_{x(t)}$ and $l_{x'(t)}$ be the half-lines in $\al L$ containing $x(t)$ and $x'(t)$ respectively. Note that $\lim_{t \to 0}  l_{x(t)} = \lim_{t \to 0}  l_{x'(t)}$ which coincides tangent cone at $0$ of $x$ and $x'$. As $t$ small enough, we may assume that the angle between $l_{x(t)}$ and $l_{x'(t)}$ is smaller than $\pi/4$. 
Let $z(t)$ and $z'(t)$ be the images of $x'(t)$ and $\varphi(x'(t))$ under the orthogonal projection onto $l_{x(t)}$. Let $\beta(t)$ denote the angle between $\overline{0 x(t)}$ and $\overline{x(t)x'(t)}$, and let $\beta'(t)$ denote the angle between $\overline{0\varphi(x)}$ and $\overline{\varphi(x(t)) \varphi(x'(t))}$ (see Figure \ref{fig:Figure1}). By reparametrizing the curves $x$ and $x'$, we may assume that  $\rho(x(t)) = \rho(x'(t)) = t$, and hence  $ \|\varphi(x(t))\|_\infty = \|\varphi(x'(t))\|_\infty = t$. By Lemma \ref{lem_val_1}, $\beta(t)$ and $\beta'(t)$ must be bounded away from $0$ when $t$ is small enough.

 \begin{figure}[h]
    \begin{tikzpicture}
    \draw (0,0) -- (6,3) node[anchor=south] {$l_{x'(t)}$}  ;
\draw (0,0) -- (7,0) node[anchor=south] {$l_{x(t)}$};
\draw[thick,dashed] (3,0) -- (3, 1.5);
\draw[thick,dashed] (4,0) -- (4,2);
\draw[red,thick] (3, 1.5)--(3.5,0);
\draw[red,thick] (4, 2)--(6,0);

\draw (0,0) node[anchor=north] {$0$};
 \draw (2.8, 1.5) node[anchor=south] {$x'(t)$};
 \fill (3, 1.5) circle (1.5pt);
  \draw (4, 2.3) node[anchor=south] {$\varphi(x'(t))$};
 \fill (4,2) circle (1.5pt);
   \draw (3.5,0) node[anchor=north] {$x(t)$};
    \fill (3.5,0) circle (1.5pt);
   \draw (6,0) node[anchor=north] {$\varphi(x(t))$};
    \fill (6,0) circle (1.5pt);
      \draw (2.7,-0.2) node[anchor=north] {$z(t)$};
       \fill (3, 0) circle (1.5pt);
      \draw (4.3,-0.2) node[anchor=north] {$z'(t)$};
       \fill (4,0) circle (1.5pt);
    \coordinate (origo) at (0,0);
    \coordinate (pivot) at (2,1);
    \coordinate (bob) at (1,0);
     \pic [draw, -, "", angle eccentricity=1.5] {angle = bob--origo--pivot};
     
     \coordinate (origo1) at (3.5,0);
    \coordinate (pivot1) at (2,0);
    \coordinate (bob1) at (3,1.5);
     \pic [draw, -, "", angle eccentricity=1.5] {angle = bob1--origo1--pivot1};
      \coordinate (origo2) at (6,0);
    \coordinate (pivot2) at (2,0);
    \coordinate (bob2) at (4,2);
     \pic [draw, -, "", angle eccentricity=1.5] {angle = bob2--origo2--pivot2};
     
     \draw (2.5,0.2) node[anchor=south] {$\beta(t)$};
      \draw (5,0) node[anchor=south] {$\beta'(t)$};
\end{tikzpicture}
\caption{}
  \label{fig:Figure1}
\end{figure}
 We have $ \|x(t)-x'(t)\|\sin \beta(t) = \|x'(t)-z(t)\|$ and $ \|\varphi(x(t)) - \varphi(x'(t))\| \sin \beta'(t) = \|\varphi(x'(t)) - z'(t)\|$. It follows that 
 $$ \frac{\|\varphi(x(t)) -\varphi(x'(t))\|}{\|x (t)- x'(t)\|} \frac{\sin \beta'(t)}{\sin \beta(t)}  = \frac{\|\varphi(x'(t)) - z'(t)\|}{\|x'(t)-z(t)\|} = \frac{\|\varphi(x'(t))\|}{\|x'(t)\|}\sim 1.$$
 
Hence, 
 $$ \frac{\|\varphi(x(t)) -\varphi(x'(t))\|}{\|x(t) - x'(t)\|} \sim  \frac{\sin \beta(t)}{\sin \beta'(t)} \sim 1$$
 since $\beta(t)$ and $\beta'(t)$ are bounded away from $0$. This again gives a contradiction.
\end{proof}


 {\it Proof of Theorem \ref{thm_main}}.\\
{\bf Case I:  $\rho(x) = \|x\|_\infty$.} \\
For $1\leq k \leq n$,  set $M^+_{k}:=\{x \in \bb R^n: \|x\|_\infty = x_k\}$ and $M^-_{k}:=\{x \in \bb R^n: \|x\|_\infty = -x_k\}$ where $x = (x_1, \ldots, x_n)$. Then $\bb R^n = \bigcup_{k=1}^n M_k^+ \cup M_k^- $. Let $\al S = \{S_i\}_{i\in I}$ be a partition of $\bb R^n$ into Lipschitz cells compatible with $\{X, \{0\}, M_1^+, M_1^-, \ldots,  M_n^+, M_n^- \}$ (see Theorem \ref{thm_ku_par}). It is clear that for  $S_i\in \al S$,  after permuting the coordinates, we have that $S_i$ is contained in the set $\al C_n$ which is defined before Lemma \ref{lem_LLNE}. It follows from Lemma \ref{lem_LLNE} that $(S_i, 0)$ is LLNE with respect to $\rho$. 

Applying Lemma \ref{lem_conic2} to $\al S$ we obtain $r_0>0, C>0$, a stratification $\Sigma$ compatible with $\al S\cap \rho^{-1}((0, r_0])$ and a stratified vector field $\xi$ on $\Sigma$ such that the flow $\Phi: $ generated by $\xi$ preserves the strata of $\Sigma$ and for $x\in X \cap \rho^{-1}((0,r_0])$, the domain of $\Phi_x$ contains the interval $[0, \rho(x))$ and for $s \in [0, \rho(x))$ we have 

     (a) $ \Phi_x (s) \in X \cap  \rho^{-1}(\rho(x)-s)$, and 
    
     (b) $\length(\Phi_x ([0, s]) \leq Cs$.
     .
     
Since $\Sigma$ is a refinement of $\al S$ and $\Phi$ preserves the strata of $\Sigma$, $\Phi$ preserves the strata of $\al S$. By putting $\Phi_x(\rho(x)): = 0$, we have $\Phi_x$ is a continuous curve passing through $x$ and ending at the origin.

Let  $p, q \in X \cap \rho^{-1} ([0,r_0])$. Set $r_1 := \|p\|_{\infty}, r_2 :=\|q\|_{\infty}$ and $s := |r_1 - r_2|$.  We may assume that $r_1 \leq r_2$.  We denote by $X(r) := X\cap \rho^{-1}(r)$. Let $\gamma$ be the integral curve of $\Phi$ through $p$ (i.e. $\gamma = \Phi_p$). We may view $\gamma$ as its image. Let $p' := \gamma \cap \rho^{-1}(r_2)$. It is clear that
\begin{equation}\label{equ_01}
  \|p - q\| \sim  \|p-q\|_{\infty}  \geq  | \|p\|_{\infty} - \|q\|_{\infty}| = s.
\end{equation}
By (b), 
\begin{equation}\label{equ_02}
      \|p -p'\| \leq d_{\gamma} (p, p')\leq Cs.
\end{equation}
Thus, 
\begin{equation}\label{equ_03}
    \|q - p'\| \lesssim \|p - q\| + \|p - p'\| \lesssim \|p - q\| .
\end{equation}

We will prove that $(iii) \Rightarrow (ii) \Rightarrow (i) \Rightarrow (iii).$
Let us first  prove that $(iii) \Rightarrow (ii)$. Given two points $p, q$ in $X$ as above. Suppose that they are both in $X_i$. We have
\begin{align*}
d_X (p, q) & \leq d_\gamma(p, p') + d_{X_i(r_2)} (p', q) \\
& \lesssim \|p - q\| + K \|q- p'\| \text{ (by (\ref{equ_01}),(\ref{equ_02}) and the fact that $X_i$ is LLNE w.r.t $\rho$)},\\
& \lesssim   \|p - q\| \text{(by (\ref{equ_03}))},
\end{align*}
where $K$ is a LLNE constant for $X_i$ at $0$. This shows that $X_i$ is LNE at $0$.

Next we prove $(ii) \Rightarrow (i)$. Because all $(X_i, 0)$  are LNE, to prove $(X,0)$ is LNE,  it suffices to check the LNE condition for two points $p \in X_i$ and $q\in X_j$, $i\neq j$. Let $p'$, $r_1$, $r_2$ be defined as above.  Let $\tilde{\gamma}$ be the integral curve of $\Phi$ through $q$.  It follows from (\ref{equ_01}), (\ref{equ_02}) and (\ref{equ_03}) that 
$$
   \|p - q\| \sim \| p - p'\| + \| p'- q\|. 
$$
Since $p'\in X_i(r_2)$ and $q\in X_j(r_2)$, by the hypothesis
$$    \|p' - q\| \geq \dist(X_i(r_2), X_j(r_2)) \geq C r_2.
$$
Therefore, 
\begin{equation}
   \|p - q\|  \gtrsim r_2. 
\end{equation}
On the other hand, 
$$d_X(p, q) \leq d_{\gamma}(p, 0) + d_{\tilde{\gamma}} (0, q) \lesssim r_1 + r_2 \lesssim r_2.$$
Thus, $\|p - q\| \gtrsim d_X(p, q)$. This implies that $X$ is LNE.

We now prove $(i) \Rightarrow (iii)$.

\underline{{\it Case I.1:}} $(X\setminus\{0\},0)$ has only one connected component. 

Let $p$ and $q$ be two points in $X(r)$ and let  $\beta$ be a curve in $X$ such that $\length(\beta) = d_X(p, q)$. Taking $r$ small enough we can assume that $\beta$ is contained in $X\cap \rho^{-1}([0,r_0])$.  

If $\beta$ goes through the origin, then $\length(\beta) \gtrsim r$. By the LNE of $(X,0)$, 
$$\length(\beta) \lesssim \|p -q\| \lesssim r$$ 
(note that for any points $x,y$ in $\bb R^n\cap \rho^{-1}(r)$ we have $\|x - y\|\lesssim r$).
This yields that 
\begin{equation}\label{equ_3.3}
    \length(\beta) \sim \|p -q\|\sim r.
\end{equation}

Let $\beta'$ be a curve in $X(r)$ such that $\length (\beta') = d_{X(r)}(p,q)$. Define $\Lambda' := \{i: \beta' \cap \overline{S}_i(r) \neq \emptyset\}$.  We denote by $p'_{i,1}, p'_{i,2} \in \overline{S}_i$  the starting and ending points of $\beta'$ in $\overline{S}_i(r)$, $i \in \Lambda'$. It follows from the LLNE of $(S_i,0)$ that there are rectifiable curves  $\alpha_i: [0,1] \to \overline{S_i}(r)$ such that $\alpha_i(0) =p'_{i,1}$, $\alpha_i(1) = p'_{i,2}$, $\alpha_i((0,1)) \subset S_i$ and $\length(\alpha_i) \sim  \|p'_{i,1}- p'_{i,2}\|$ (see (2), Remark \ref{rem_L_cell}). Let $\alpha$ be the union of  $\alpha_i$'s, we then have
   $$ d_{X(r)}(p,q) \leq \length(\alpha) = \sum_{i\in \Lambda'} \length(\alpha_i)
    \lesssim \sum_{i\in \Lambda'}\|p'_{i,1}-p'_{i,2}\|\lesssim r.$$
Combining with (\ref{equ_3.3}) we see that $X$ is LLNE with respect to $\rho$.

We now consider the case that the curve $\beta$ does not go through the origin. Set $\Lambda := \{i: \beta \cap \overline{S}_i \neq \emptyset\}$.  We denote by $p_{i,1}, p_{i,2} \in \overline{S}_i$ respectively the starting  and the ending points  of $\beta$ in $\overline{S_i}$, $i\in \Lambda$. Let $\beta_{i, j}$, $i \in \Lambda$, $j = 1, 2$ be the integral curves through $p_{i, j}$  generated by $\xi$ and let  $q_{i, j} := \beta_{i, j} \cap X (r)$ (see Figure \ref{fig:Figure2}).
 \begin{figure}[h]
\begin{tikzpicture}
    \draw[->] (0,0) -- (5,0) node[anchor=north] {$x_1$};
    \draw[->] (0,0) -- (0,5) node[anchor=east] {$\mathbb{R}^{n-1}$};
    \draw (0,0) node[anchor=north] {$0$};
    \draw[thick] plot [smooth,tension=1] coordinates{ (0,0) (0.2,1) (0.5,4) };
    \draw[thick] plot [smooth,tension=1] coordinates{ (0,0) (2,0.5) (4,0.5) };
\draw[fill=gray, gray] (0,0) -- (4,2.5) -- (4,4) -- cycle;
\draw[] (0,0)-- (4,2.5);
\draw[] (0,0)-- (4,4);

\draw[dashed, thick] (0,4)-- (4,4)--(4,0);

\draw[thick, blue] plot [smooth,tension=1] coordinates{ (1,4) (1.4,1.7) (4,1) };
\draw (1,4) node[anchor=south] {$p$};
\draw (4,1) node[anchor=west] {$q$};
\draw (4,1) node[anchor=west] {$q$};
\draw (4,4) node[anchor=south] {$q_{i, 1}$};
\draw (4,2.5) node[anchor=west] {$q_{i, 2}$};
\draw (2,2)  node[anchor=east] {$p_{i, 1}$};
\draw (1.9,1.17) node[anchor=north] {$p_{i, 2}$};
\draw (1,3) node[anchor=east] {$\beta$};
\draw (3,3) node[anchor=north] {$S_i$};
\fill (4,1) circle (1.5pt);
\fill (1,4) circle (1.5pt);
\fill (2.05,1.3) circle (1.5pt);
\fill (1.56,1.56) circle (1.5pt);
\fill (4,2.5) circle (1.5pt);
\fill (4,4) circle (1.5pt);

\draw (4,0) node[anchor=north] {$r$};
\draw[->] (5,2) --(4,2);
\draw (5,2) node[anchor=west] {$\rho^{-1}(r)$};
\end{tikzpicture}
\caption{}
  \label{fig:Figure2}
\end{figure}

Set $ l := \max_{i\in \Lambda} \{ \|p_{i,j} - q_{i,j}\|\}$. We may assume $l = \|p_{i,1} - q_{i,1}\|$ for some $i \in \Lambda$. By (b), we have 
\begin{equation*}\label{equ3}
    l \leq C |\|p_{i,1}\|_{\infty} - \|q_{i,1}\|_{\infty}|,
\end{equation*}
where $C$ is a universal constant  depending only on $X$.

On the other hand, 
\begin{equation*}\label{equ4}
    \length(\beta) \geq d_{\beta} (p_{i,1}, p) \geq \| p_{i,1} - p\| \sim \|p_{i,1} - p\|_{\infty} \geq |\|p_{i,1}\|_{\infty} - \|p\|_{\infty}|.
\end{equation*}
Since $p$ and  $q_{i,1}$ lie in $\rho^{-1}(r)$, $\|p\|_{\infty} = \|q_{i,1}\|_{\infty}$. It follows that
\begin{equation}\label{equ5}
    \length(\beta) \geq |\|p_{i,1}\|_{\infty} - \|q_{i,1}\|_{\infty}|\gtrsim l.
\end{equation}
Since $X$ is LNE at $0$, we have  
\begin{equation}\label{equ6}
   l \lesssim \length(\beta) \sim \|p - q\|.
\end{equation}
In addition, 
\begin{equation}\label{equ7}
    \|q_{i,1} - q_{i,2}\| \leq \|q_{i, 1} - p_{i,1}\| + \| p_{i, 1} - p_{i, 2}\|  + \| p_{i,2}- q_{i, 2}\|
     \leq \| p_{i, 1} - p_{i, 2}\| + 2 l.
\end{equation}
Since $\overline{S_i}$ is LNE,  
\begin{equation}\label{equ8}
    \length(\beta) = \sum_{i\in \Lambda} d_{\overline{S}_i} (p_{i,1}, p_{i, 2}) \sim \sum_{i\in \Lambda}  \|p_{i,1} -p_{i, 2}\|.
\end{equation} 
 From (\ref{equ6}), (\ref{equ7}) and  (\ref{equ8}) and the fact that $(S_i,0)$ is LLNE with respect to $\rho$, we get
\begin{align*}
    d_{X(r)} (p, q)  &\lesssim \sum_{i\in \Lambda} d_{\overline{S}_{i}(r)} (q_{i,1}, q_{i, 2})  \sim  \sum_{i\in \Lambda}  \|q_{i,1} -q_{i, 2}\| \\
    &\leq  \sum_{i\in \Lambda}  (\|p_{i,1} -p_{i, 2} \|+ 2l)  \lesssim (2m+1) \length(\beta) \sim (2m+1)\|p-q\|, 
\end{align*} 
where $m := \# \Lambda$.
This implies that $X$ is LLNE with respect to $\rho$. 

\underline{{\it Case I.2:}} $(X\setminus\{0\},0)$ has more than one connected component. 

By Case 1, $X_i$ are LLNE with respect to $\rho$. Thus,  the first condition in $(iii)$ is satisfied.  We assume on the contrary that the second condition in $(iii)$ fails i.e., there are $j, k \in I$, $j\neq k$ such that 
\begin{equation}\label{equ7.1}
    \dist(X_j(r), X_k(r))\ll r.
\end{equation}
Let $\tilde{\gamma}_j: (0, \varepsilon) \to X_j$ and $\tilde{\gamma}_k: (0, \varepsilon)\to X_k$  be  definable curves in $X_j$ and $X_k$ respectively such that $\lim_{r\to 0}\tilde{\gamma}_j(r) = \lim_{r\to 0}\tilde{\gamma}_k(r) =0$, $\rho(\tilde{\gamma}_j(r)) =  \rho( \tilde{\gamma}_k(r)) = r$ and $\|\tilde{\gamma}_j(r)  - \tilde{\gamma}_k(r)\| \sim \dist(X_j(r), X_k(r))$. Note that every rectifiable curve connecting $\tilde{\gamma}_j(r)$ and $ \tilde{\gamma}_k(r)$ must go through the origin, and hence its length is $\gtrsim r$. By (\ref{equ7.1}), we have
$$\lim_{r\to 0} \frac{d_X(\tilde{\gamma}_j(r), \tilde{\gamma}_k(r))}{ \|\tilde{\gamma}_j(r) - \tilde{\gamma}_k(r)\|} \to \infty.$$
This contradicts the fact that $X$ is LNE.

{\bf Case II:  $\rho$ is general.}

First, we claim that it is possible to extend $\rho$ to a definable radius function on a neighbourhood of the origin in $\bb R^n$. Indeed, since $\rho(x) \sim \|x\|$, there is $K >0$  such that $\frac{1}{K} \|x\| \leq \rho(x) \leq K\|x\|$ for every $x \in X$ near the origin. Let $\tilde{\rho}$ be a definable Lipschitz extension of $\rho$ to the whole of $\bb R^n$. It is easy to check that $\overline{\rho}(x) := \max\{\frac{1}{K}\|x\|, \min \{ K\|x\|, \tilde{\rho}\}\}$ is the desired extension. 

Let us  use the same notation $\rho$ for such an extension.  Let $\varphi: (\bb R^n, 0) \to (\bb R, 0)$ defined by $\varphi(x) := \frac{\rho(x)}{\|x\|_{\infty}} x$ if $x \neq 0$ and $\varphi(0): = 0$. It is obvious that $\varphi (L_{\rho, r}(X)) = L_{\rho', r}(\varphi(X))$ where $\rho'(x) = \|x\|_{\infty}$. Moreover, by Lemma \ref{lem_bi-Lipschitz}, $\varphi$ is a germ of a  bi-Lipschitz definable homeomorphism. Then, we have the following equivalent statements: 

(1) $X$ is LNE.

(2) $\varphi(X)$ is LNE.

(3) $(\varphi(X_i), 0)$ is LNE for every $i \in I$ and there is  $C>0$ such that $$\dist(L_{\rho', r}(\varphi(X_j)),L_{\rho', r}(\varphi(X_k)) \geq C r$$ for all $j, k \in I$, $j\neq k$ and $r$ sufficiently small.

(4) $(\varphi(X_i),0)$ is LLNE with respect to $\rho'$ for ever $i \in I$ and there is  $C>0$ such that $$\dist(L_{\rho', r}(\varphi(X_j)),L_{\rho', r}(\varphi(X_k)) \geq C r$$ for all $j, k \in I$, $j\neq k$ and $r$ sufficiently small.

(5) $(X_i, 0)$ is LNE  for every $i \in I$ and there is  $C>0$ such that $$\dist(L_{\rho, r}(X_j),L_{\rho, r}(X_k) \geq C r$$ for all $j, k \in I$, $j\neq k$ and $r$ sufficiently small.

(6) $(X_i,0)$ is LLNE with respect to $\rho$ for every $i \in I$ and there is  $C>0$ such that $$\dist(L_{\rho, r}(X_j),L_{\rho, r}(X_k) \geq C r$$ for all $j, k \in I$, $j\neq k$ and $r$ sufficiently small.

Indeed, $(1) \Leftrightarrow (2)$ follows from the definition of LNE;  and 
$(2) \Leftrightarrow (3) \Leftrightarrow (4)$ follows from Case I; 
$(3) \Leftrightarrow (5)$  and $(4) \Leftrightarrow (6)$ follow from the fact that $\varphi$ brings $X_i$ to $\varphi(X_i)$ (as germs at $0$) and  $L_{\rho, r} (X)$ to $L_{\rho', r} (\varphi(X))$. This completes the proof. 

\begin{rem}\label{rem_final} \rm 
The Lipschitzness of the function $\rho$ in Theorem \ref{thm_main} is necessary even for the semialgebraic case. For example, let  $X: = \{(x_1, x_2) \in \bb R^2: x_2^2 \leq  \frac{x_1^3}{2}\}$ and $\rho: X \to \bb R$ be a function germ defined by $\rho(x_1, x_2) := x_1 - \sqrt[3]{x_2^2}$. It is obvious that  $\rho$ is a continuous semialgebraic satsfying $\rho(x) \sim \|x\|$ over $(X,0)$. The function $\rho$ is not Lipschitz since for $z_1 = (x_1, x_1^{3/2})$ and $z_2 = (x_1, 0)$, we have $$\|\rho(z_1) - \rho(z_2)\| = |x_1|\gg x_1^{3/2} = \|z_1 - z_2\|.$$
The $r$-link of $(X, 0)$ associated to $\rho$: $L_{\rho, r}(X) = \{(x, y) \in X: y^2 = (x - r)^3\}$ is clearly not LNE.

Note that it is always possible to extend $\rho$ to a continuous semialgebraic function $\rho':(\bb R^n, 0) \to \bb R$ such that $\rho(x) \sim \|x\|$. To do it, one first extends $\rho$ to a continuous semialgebraic function on the whole of $(\bb R^n, 0)$ (this is possible since $(X, 0)$ is closed, see for example \cite{Dries}, Corollary 3.10 or \cite{Asc-Fis}, Lemma 6.6) then use the same arguments as in the extension of radius functions (see in Case II of the proof of Theorem \ref{thm_main}). 
\end{rem}

\section{Application}\label{section3}
This section presents  a counterexample to a question asked in \cite{Bobadilla} about a sufficient condition on MD-homologies for a subanalytic germ being  Lipschitz normally embedded (see Question \ref{question}). We refer the reader to \cite{Bobadilla} for  the precise definition of MD-homology of a subanalytic germ and its basic properties. Here, we just recall some necessary facts that help to understand the question.

Let  $f: (X, x_0, d_1) \to (Y, y_0, d_2)$ be a subanalytic map between two metric subanalytic germs. We call $f$ {\it linearly vertex approaching} (or {\it l.v.a}) if there is a constant $K>0$ such that 
$$ \frac{1}{K} d_1(x, x_0) \leq d_2(f(x), y_0) \leq K d_1(x, x_0).$$

If $f$ is Lipschitz l.v.a then it induces a homomorphism $f_*: MDH_{\bullet}^b ((X, x_0, d_1), A)  \to MDH_{\bullet}^b ((Y, y_0, d_2), A)$ between $b$-MD-homologies of $(X, x_0)$ and $(Y, y_0)$ for every $b \in (0, +\infty]$ where $A$ is an abelian group. In particular, if $f$ is a bi-Lipschitz subanalytic homeomorphism then $f_*$ is an isomorphism.  

We consider the identity map $\id_{(X,x_0)}: (X, x_0, d_{inn}) \to (X, x_0, d_{out})$. It is clear that $\id_{(X, x_0)}$ is a Lipschitz l.v.a subanalytic map.  If  $(X, x_0)$ is LNE then $\id_{(X, x_0)}$ is a bi-Lipschitz subanalytic homeomorphism, therefore
$$\id_{(X, x_0) *}: MDH_{\bullet}^b ((X, x_0, d_{inn}), A)  \to MDH_{\bullet}^b ((X, x_0, d_{out}), A)$$
is an isomorphism for every $b \in (0, \infty]$. The question is whether the converse holds. More precisely that:

\begin{question}[\cite{Bobadilla}, Problem 147 ]\label{question}  
  Suppose that for every $x\in (X,x_0)$, the map $$\id_{(X,x)*}: MDH_{\bullet}^b ((X, x, d_{inn}), A)  \to MDH_{\bullet}^b ((X, x, d_{out}), A)$$ is an isomorphism. Is $(X, x_0)$ LNE ?
\end{question}

The following example shows that in general the answer is negative. 

\begin{example}
Let $X:=\{(t, x, z) \in \bb R^3, z^2 = t^2x^2, 0\leq x \leq t\}$. Then 
\begin{enumerate}
    \item $(X, 0)$ is not LNE.
    \item $\id_{(X,x)*}: MDH^b_\bullet (X, x, d_{inn}) \to MDH^b_\bullet (X, x, d_{out})$ is an isomorphism for every $x \in (X,0)$.
\end{enumerate}
Therefore, Question \ref{question} has a negative answer.
\end{example}

\begin{proof}
Let $\rho (t, x, z)  = t$. It is clear that the restriction of $\rho$ to $(X, 0)$ is a radius function. Set $X(r) := X\cap \{\rho^{-1}(r)\}$. It is easy to see that $X(r)$ is LNE with the LNE constant $\sim 1/r$, which tends to $\infty$ when $r$ tends to $0$. This means  $(X,0)$ is not LLNE with respect to $\rho$. Note that $(X,0)$ is connected, so by Theorem \ref{thm_main_2} (see also Theorem \ref{thm_Men-Sam_2}) it is not LNE. Thus, $(1)$ is proved.

We now show (2). Let $S_0 := \{0\}$, $S_1 := \{(t,0,0) \in \bb R^3, t>0\}$, $S_2:=\{(t, x, z)\in X, 0<x = t\}$ and $S_3 := X \setminus \cup_{i=0}^2 S_i$ (see Figure \ref{fig:Figure3}). 

\begin{figure}[h]
\begin{tikzpicture}
    \draw[->] (0,0) -- (5,0) node[anchor=north] {$t$};
    \draw[->] (0,0) -- (0,3) node[anchor=east] {$x$};
    \draw[->] (0,0) -- (-1.5,-1) node[anchor=east] {$z$};
    \draw (0,0) node[anchor=north] {$0$};
    
     \draw[thick, dashed] (4,0) -- (3.7,1);
     \draw[thick, dashed] (4,0) -- (4.2,1.7);
     \draw[thick] (0,0) -- (4,0);

    \draw[thick, blue ] plot [smooth,tension=1] coordinates{ (0,0) (2,0.3) (3.7,1) };
    \draw[ thick, blue] plot [smooth,tension=1] coordinates{ (0,0) (2,0.5) (4.2,1.7) };

    \draw (2,2) node[anchor=south] {$X$};
    \draw (4,0) node[anchor=north] {$r$};
    
      \draw[->] (3,-0.7) --(3,0);
       \draw (3,-0.7) node[anchor=north] {$S_1$};
       
      \draw[->] (3.5,2.5) --(3.5, 1.3);
       \draw[->] (3.5,2.5) --(3.3, 0.8);
       \draw (3.5,2.5) node[anchor=south] {$S_2$};
    
    \fill (0,0) circle (1.5pt);
    
    \draw[->] (7,0) -- (10,0) node[anchor=north] {$z$};
    \draw[->] (8.5,-1) -- (8.5,3) node[anchor=east] {$x$};
    \draw[thick] (8,1.5)--(8.5,0)--(9,1.5);
    \draw[dashed](6.8, 3.3)--(10.2, 3.3)--(10.2, -1.2)--(6.8, -1.2) --(6.8, 3.3);
     \draw (8.5,-2) node[anchor=north] {$X\cap \rho^{-1}(r)= \{x = \frac{1}{r}|z|\}$};
    
\end{tikzpicture}
\caption{}
  \label{fig:Figure3}
\end{figure}

It is clear that $\{S_i\}_{i=0}^3$ is a stratification of $X$. Let $x \in X$. If $x \in S_3$, then $x$ is a smooth point, hence $(X, x)$ is LNE, hence (2) is satisfied. If $x \in S_2$, the germ $(X,x)$ is a smooth manifold with boundary, so it is also LNE hence $(2)$ is again true. If $x \in S_1$, by Valette's Lipschitz Triviality Theorem (see \cite{Valette1}, Theorem 2.2) in a neighbourhood of $x$, $X$ is bi-Lipschitz equivalent to $X(r)\times (0, \varepsilon)$ where $r = \|x\|$, which is obviously LNE. This implies that $(2)$ holds.  The only case that needs verifying is  $x=0$. 

Consider the following map: 
$$H: X \times I \to X, (t, x, z, s) \mapsto H_s(t, x,z) := (t, sx, sz).$$  
We show  that $H$ is a Lipschitz l.v.a metric homotopy for both the outer  and the inner metrics. Fix $s$ and let  $w = (t_1, x_1, z_1)$ and $w'=(t_2, x_2, z_2)$. We have 
\begin{align*}
    \|H_s(w) - H_s(w')\|  & = \|(t_1, sx_1, sz_1) - (t_2, sx_2, sz_2)\|  \\ 
    & \leq |t_1 - t_2| + s(|x_1 - x_2| + |z_1 - z_2|)\lesssim \|w - w'\|.
\end{align*}
This shows that $H_s$ is Lipschitz with respect to the outer metric.

Observe that  $X$ consists of two branches $X_1 := \{(t, x, z)\in X, z\geq 0\}$ and $X_2 :=\{(t, x, z)\in X, z\leq 0\}$ and each branch is LNE. Moreover, $H_s(-)$ preserves these branches. 

If $w$ and $w'$ are in the same branch we have 
$$d_{inn} (H_s(w), H_s(w')) \sim \|H_s(w) - H_s(w')\| \lesssim \|w - w'\| \sim d_{inn}(w, w').$$
Now assume that $w\in X_1$ and $w'\in X_2$. We have
\begin{align*}
    d_{inn} (H_s(w), H_s(w')) &= d_{inn} ((t_1, sx_1, sz_1), (t_2, sx_2, sz_2))\\
    & \leq  d_{inn} ((t_1, sx_1, sz_2), (t_2, 0, 0)) + d_{inn} ((t_2, 0, 0), (t_2, sx_2, sz_2))\\
    &\lesssim |t_1 - t_2| + s(|x_1| + |z_1|) + s(|x_2| + |z_2|).
\end{align*}
Let $\gamma$ be a curve connecting $w$ and $w'$ which realizes the inner distance between $w$ and $w'$. Since $w$ and $w'$ lie in two different branches of $X$, $\gamma$ has to pass through the $t$-axis. Hence $$d_{inn}(w, w') = \length(\gamma) \gtrsim |t_1 - t_2| + (|x_1| + |z_1| + |x_2|+|z_2|).$$
It follows that $ d_{inn} (H_s(w), H_s(w')) \lesssim d_{inn}(w, w')$. Consequently,  $H_s$ is Lipschitz with respect to the inner metric. Note that Lipschitz constants of $H_s$ can be chosen to be independent of $s$.

Let $Y: = \{(t, 0, 0), t \geq 0\}$ and $g: X \to Y, w \mapsto g(w) = H_0(w)$. Clearly, $H_0 = \iota\circ g$ where $\iota: Y \to X$ is the inclusion map. Since $H_1= \id_{(X,0)}$ and $H_0$ are metrically homotopic for both the inner  and the outer metrics,  $\id_{(X,0) *},  H_{0 *}:  MDH^b_{\bullet} ((X, 0, d), A) \to MDH^b_{\bullet} ((X, 0, d), A)$ (where $d \in \{d_{inn}, d_{out} \}$) represent the same homomorphism which is actually an isomorphism (see \cite{Bobadilla},  Theorem 81 (2)). Since $H_{0 *} = \iota_*\circ g_* $, which is an isomorphism, and $i_*$ is injective,  $g_*$ must be an isomorphism.  We have the following commutative diagram: 

\[ \begin{tikzcd}
MDH^b_{\bullet}((X, 0, d_{inn}), A) \arrow{r}{g_*} \arrow[swap]{d}{\id_{(X,0) *}} & MDH^b_{\bullet}((Y, 0, d_{inn}), A) \arrow{d}{\id_{(Y,0)*}} \\%
MDH^b_{\bullet}((X, 0, d_{out}), A)\arrow{r}{g_*}& MDH^b_{\bullet}((Y, 0, d_{out}), A)
\end{tikzcd}
\]

Since $(Y,0)$ is a germ of a half line which is LNE, $\id_{(Y,0) *}$ is an isomorphism. Thus,  $\id_{(X,0) *}$ is also an isomorphism. 
\end{proof}

\begin{rem}\rm
Since our counterexample is a set with non-isolated singularities, Question \ref{question} is still open in  the isolated singularities case.
\end{rem}

\subsection*{Acknowledgements} We would like to thank the referee for valuable remarks and suggestions. We would also like to thank Javier Bobadilla for his interest and useful discussion on the example. The research is funded  by Vietnam National Foundation for Science and Technology Development (NAFOSTED) under the grant number 101.04-2019.316. It was also  supported by the ERCEA 615655 NMST Consolidator Grant,  by the Basque Government through the BERC 2018--2021 program and by the Spanish Ministry of Science, Innovation and Universities: BCAM Severo Ochoa accreditation SEV-2017--0718.

\bibliographystyle{siam}
\bibliography{Biblio}

\begin{thebibliography}{10}

\bibitem{Asc-Fis}
{\sc M.~Aschenbrenner and A.~Fischer}, {\em Definable versions of theorems by
  {K}irszbraun and {H}elly}, Proc. Lond. Math. Soc. (3), 102 (2011),
  pp.~468--502.

\bibitem{Anne2}
{\sc A.~Belotto~da Silva, L.~Fantini, and A.~Pichon}, {\em On lipschitz
  normally embedded complex surface germs}, Preprint,
  https://arxiv.org/abs/2006.01773,  (2020).

\bibitem{Lev2}
{\sc L.~Birbrair, A.~Fernandes, and W.~D. Neumann}, {\em On normal embedding of
  complex algebraic surfaces}, in Real and complex singularities, vol.~380 of
  London Math. Soc. Lecture Note Ser., Cambridge Univ. Press, Cambridge, 2010,
  pp.~17--22.

\bibitem{Bir-Men}
{\sc L.~Birbrair and R.~Mendes}, {\em Arc criterion of normal embedding}, in
  Singularities and foliations. geometry, topology and applications, vol.~222
  of Springer Proc. Math. Stat., Springer, Cham, 2018, pp.~549--553.

\bibitem{Bir-Men-Nun}
{\sc L.~Birbrair, R.~Mendes, and J.~J. Nu\~{n}o Ballesteros}, {\em Metrically
  un-knotted corank 1 singularities of surfaces in {$\Bbb{R}^4$}}, J. Geom.
  Anal., 28 (2018), pp.~3708--3717.

\bibitem{Bir-Mos}
{\sc L.~Birbrair and T.~Mostowski}, {\em Normal embeddings of semialgebraic
  sets}, Michigan Math. J., 47 (2000), pp.~125--132.

\bibitem{Bobadilla}
{\sc J.~F. Bobadilla, S.~Heinze, M.~P. Pereira, and J.~E. Sampaio}, {\em
  Moderately discontinuous homology}, To appear in Communications on pure and
  applied Mathematics, https://arxiv.org/abs/1910.12552.

\bibitem{Coste}
{\sc M.~Coste}, {\em An introduction to o-minimal geometry}, Dip. Mat. Univ.
  Pisa, Dottorato di Ricerca inMatematica, Istituti Editoriali e Poligrafici
  Internazionali, Pisa, 2000.

\bibitem{Den-Tib}
{\sc M.~Denkowski and M.~Tib\u{a}r}, {\em Testing {L}ipschitz non-normally
  embedded complex spaces}, Bull. Math. Soc. Sci. Math. Roumanie (N.S.),
  62(110) (2019), pp.~93--100.

\bibitem{Fer-Sam1}
{\sc A.~Fernandes and J.~E. Sampaio}, {\em Tangent cones of {L}ipschitz
  normally embedded sets are {L}ipschitz normally embedded. {A}ppendix by
  {A}nne {P}ichon and {W}alter {D}. {N}eumann}, Int. Math. Res. Not. IMRN,
  (2019), pp.~4880--4897.

\bibitem{Fer-Sam2}
\leavevmode\vrule height 2pt depth -1.6pt width 23pt, {\em On {L}ipschitz
  rigidity of complex analytic sets}, J. Geom. Anal., 30 (2020), pp.~706--718.

\bibitem{Ka}
{\sc K.~Katz, M.~Katz, D.~Kerner, and Y.~Liokumovich}, {\em Determinantal
  variety and normal embedding}, J. Topol. Anal., 10 (2018), pp.~27--34.

\bibitem{Kur}
{\sc K.~Kurdyka}, {\em On a subanalytic stratification satisfying a {W}hitney
  property with exponent {$1$}}, in Real algebraic geometry ({R}ennes, 1991),
  vol.~1524 of Lecture Notes in Math., Springer, Berlin, 1992, pp.~316--322.

\bibitem{Kur-Par}
{\sc K.~Kurdyka and A.~Parusi\'{n}ski}, {\em Quasi-convex decomposition in
  o-minimal structures. {A}pplication to the gradient conjecture}, in
  Singularity theory and its applications, vol.~43 of Adv. Stud. Pure Math.,
  Math. Soc. Japan, Tokyo, 2006, pp.~137--177.

\bibitem{Loi2}
{\sc T.~L. Loi}, {\em Whitney stratification of sets definable in the structure
  {$\bold R_{\exp}$}}, in Singularities and differential equations ({W}arsaw,
  1993), vol.~33 of Banach Center Publ., Polish Acad. Sci. Inst. Math., Warsaw,
  1996, pp.~401--409.

\bibitem{Loi1}
\leavevmode\vrule height 2pt depth -1.6pt width 23pt, {\em Verdier and strict
  {T}hom stratifications in o-minimal structures}, Illinois J. Math., 42
  (1998), pp.~347--356.

\bibitem{Men-Sam}
{\sc R.~Mendes and J.~E. Sampaio}, {\em On {L}ink of {L}ipschitz normally
  embedded sets}, https://arxiv.org/pdf/2101.05572.pdf,  (2021).

\bibitem{Anne1}
{\sc F.~Misev and A.~Pichon}, {\em Lipschitz normal embedding among
  superisolated singularities}, International Mathematics Research Notices,
  2021 (2021), pp.~13546--13569.

\bibitem{Neu-Per-Pic1}
{\sc W.~D. Neumann, H.~M.~l. Pedersen, and A.~Pichon}, {\em A characterization
  of {L}ipschitz normally embedded surface singularities}, J. Lond. Math. Soc.
  (2), 101 (2020), pp.~612--640.

\bibitem{Neu-Per-Pic2}
\leavevmode\vrule height 2pt depth -1.6pt width 23pt, {\em Minimal surface
  singularities are {L}ipschitz normally embedded}, J. Lond. Math. Soc. (2),
  101 (2020), pp.~641--658.

\bibitem{Ngu-Tri-Tro}
{\sc N.~Nguyen, S.~Trivedi, and D.~Trotman}, {\em A geometric proof of the
  existence of definable {W}hitney stratifications}, Illinois J. Math., 58
  (2014), pp.~381--389.

\bibitem{ng-va}
{\sc N.~Nguyen and G.~Valette}, {\em Whitney stratifications and the continuity
  of local {L}ipschitz-{K}illing curvatures}, Ann. Inst. Fourier (Grenoble), 68
  (2018), pp.~2253--2276.

\bibitem{Valette3}
{\sc G.~Valette}, {\em On subanalytic geometry, survey},
  http://www2.im.uj.edu.pl/gkw/sub.pdf.

\bibitem{Valette1}
{\sc G.~Valette}, {\em Lipschitz triangulations}, Illinois J. Math., 49 (2005),
  pp.~953--979.

\bibitem{Valette2}
\leavevmode\vrule height 2pt depth -1.6pt width 23pt, {\em The link of the germ
  of a semi-algebraic metric space}, Proc. Amer. Math. Soc., 135 (2007),
  pp.~3083--3090.

\bibitem{Dries}
{\sc L.~van~den Dries}, {\em Tame topology and o-minimal structures}, vol.~248
  of London Mathematical Society Lecture Note Series, Cambridge University
  Press, Cambridge, 1998.

\end{thebibliography}

\end{document}